\documentclass[10pt]{amsart}
\usepackage{hyperref}
\usepackage{amssymb,amsmath, amsthm, bbm}
\usepackage[all]{xy}
\usepackage{graphicx}

\theoremstyle{definition}

\newtheorem{thm}{Theorem}[section]

\newcommand{\PP}{\mathcal{P}}  
\newcommand{\ed}{\mathbf{d}}   
\newcommand{\J}{\mathbf{J}}    
\newcommand{\diag}{\mathrm{diag}}
\newcommand{\R}{\mathbb{R}}    
\newcommand{\Z}{\mathbb{Z}} 
\newcommand{\Ad}{\mathrm{Ad}}  
\newcommand{\Sl}{\mathbf{S}}   
\newcommand{\g}{\mathfrak{g}}
\newcommand{\q}{\mathfrak{q}}
\newcommand{\m}{\mathfrak{m}}
\newcommand{\be}{\begin{equation}}
\newcommand{\ee}{\end{equation}}
\newcommand{\bea}{\begin{eqnarray}}
\newcommand{\eea}{\end{eqnarray}}
\newcommand{\Proj}{\mathbb{P}} 
\newcommand{\restr}[1]{\vrule height3ex width.4pt depth1.4ex\lower1.4ex\hbox{\scriptsize $\,#1$}}
\newcommand{\rrestr}[1]{\vrule height2ex width.4pt depth0.9ex\lower0.9ex\hbox{\scriptsize $\,#1$}}
\newcommand{\e}{\mathbf{e}} 
 
\newcommand{\Null}{\mathcal{N}}

\newcommand{\zt}{\widetilde{z}}

\begin{document}

\title[]{Bifurcations of MacLaurin spheroids. A Hamiltonian Perspective}
\author{Miguel Rodr\'{\i}guez-Olmos}\thanks{Departamento de Matem\'aticas. Universitat Polit\`ecnica de Catalunya (UPC)}\thanks{miguel.rodriguez.olmos@upc.edu}

\begin{abstract} Dirichlet's problem for the dynamics of fluid bodies with ellipsoidal shape can be formulated as a Hamiltonian system invariant under the action of a symmetry Lie group.
I apply methods from Hamiltonian bifurcation theory to the study of the branch of solutions known as MacLaurin spheroids. I show that  all its bifurcations are into three types named I, $S$ and adjoint $S$ ellipsoids in agreement with previous necessary conditions obtained by Chandrasekhar by linearizing the hydrodynamic equations. 
\end{abstract}

\maketitle


\section{Introduction}

This article studies the bifurcations from the MacLaurin family of Riemann ellipsoids  using geometric methods from the theory of Hamiltonian relative equilibria. Riemann ellipsoids can be thought as the relative equilibria for a finite dimensional symmetric Hamiltonian system whose Hamilton's equations correspond to Lebovitz's formulation of Dirichlet's problem. This problem consists in studying the dynamics of a self gravitating incompressible fluid mass with constant density that must remain an ellipsoid at all times and for which the velocities are linear functions of the coordinates. (see \cite{Cha87} and \cite{RoSD99}). 

In this setting I will consider the MacLaurin family of relative equilibria solutions which consists of a branch of oblate spheroidal configurations of the fluid  where both angular momentum and circulation are parallel to the symmetry axis and the same is true for the angular velocity and internal vorticity. The elements in the MacLaurin family can be parametrized via the eccentricity of its spheroidal configuration. Different MacLaurin spheroids with the same eccentricity are regarded as the same relative equilibrium, since under the projection by the symmetry group of the system they are mapped to a single fixed point of the reduced dynamics.

The limit cases of eccentricity 0 and 1 correspond to the spherical stable equilibrium and the idealized two-dimensional disc. The spherical solution is stationary and, in an inertial frame, MacLaurin spheroids with non-zero eccentricity can be seen as rotating spheroids around the symmetry axis with  a frequency that depends on the difference of angular velocity and vorticity. As the eccentricity increases from zero, these MacLaurin spheroids are stable, and once passed the critical value $e_\mathrm{crit}\simeq 0.952887$ they become unstable. Bifurcations occur for all values of the eccentricity, but differ if the eccentricity is above, below or at this critical value.

In this work I study the problem of the existence of Riemann ellipsoids in a neighborhood of the MacLaurin family and how they are organized in branches that bifurcate from it. Dirichlet's problem has been extensively studied from many angles, including a large body of classic work within the field of analytical dynamics, and all the different types of Riemann ellipsoids as well as many of their dynamical properties are well known. A modern comprehensive account of the available results until the 1970's is given by Chandrasekhar in \cite{Cha87}. There the author obtains all Riemann ellipsoids that are solutions of Dirichlet's problem among other results that go beyond the ellipsoidal restriction, like the pear-shaped solutions discovered by Poincar\'e (see \cite{Poi85}). It also studies their linear stability. This stability analysis  has since then been advanced in several ways: In \cite{FaDe01} the authors study the stronger notion of Nekhoroshev stabiliy, which is a kind of nonlinear stability for finite but very long times. A study of the nonlinear stability, specifically addressing the $G_\mu$-stability introduced in \cite{Pat92}, can be found in \cite{ROSD07} but only for the particular solutions where the configuration ellipsoids have at least one axis of symmetry. In \cite{Cha87} the relative equilibrium solutions are organized into parametrized branches, the MacLaurin family being one of them. A linearization analysis is employed to identify points satistying a necessary condition for a bifurcation into a different family as the parameter reaches a certain value. The conditions obtained are not sufficient, and the existence of the bifurcation is later established by explicitely solving the system's differential equations. A recent normal form approach can be seen in \cite{FaPaYa23} where the authors study nonlinear bifurcations for those solutions which have configurations with 3 different axes (although their study includes one family of axisymmetric solutions, which has been referred to as `prolate spheroids' in \cite{ROSD07} and 48.(e) of \cite{Cha87}). In particular \cite{FaPaYa23} does not consider the bifurcations  from the MacLaurin family or the spherical equilibrium. In this article I address the first of these two families via a nonlinear bifurcation analysis from the point of view of geometric mechanics and Hamiltonian systems invariant under the action of a symmetry group. I will use the results of \cite{RO20, ROSD07} in order to obtain necessary and sufficient conditions for the existence of bifurcations from the MacLaurin familiy into the families of Riemann ellipsoids known as type $S$, type adjoint $S$ and type $I$, to be introduced later. A future work will address the problem of the bifurcations from the spherical equilibrium. In the remainder of this section I introduce the formal details about Dirichlet's problem and symmetric Hamiltonian systems needed in the article following the notation of \cite{ROSD07}, \cite{ROSD02} and \cite{RoSD99}.

Dirichlet's model can be seen, and this is a profound observation due to Lebovitz, as a dynamical system on the configuration space $\mathrm{SL}(3)$. We will identify it here with the group of $3\times 3$ matrices of determinant 1. The kinetic energy is obtained from the metric given by the trace 
\begin{equation}\label{metricsl3}\ll A, B \gg = T\mathrm{tr}(A^TB),\quad A,B\in T_F\mathrm{SL}(3)\end{equation}
where  $T=\frac{4\pi}{15}\rho_0$ and $\rho_0$ is the fluid's constant density. This induces a diffeomorphism of $T\mathrm{SL}(3)$ with $T^*\mathrm{SL}(3)$ and a fiberwise inner product on the latter.
The potential energy function is given by 
    \begin{equation}\label{gravpot}V(F)=-R\int_0^\infty \frac{ds}{\sqrt{s^3+I_1(F)s^2+I_2(F)s+1}}\end{equation}
where $R=\frac{8}{15}\pi^2 G\rho_0^2$, $G$ is the gravitational constant and

\begin{eqnarray*} 
    I_1(F) & = & \mathrm{tr}(FF^T)\\
    I_2(F) & = &\frac 12 (\mathrm{tr}^2(FF^T)-\mathrm{tr}(FF^TFF^T))
\end{eqnarray*}

The Lie group of symmetries in Dirichlet's problem is given by the semidirect product $G = \mathbb{Z}_2^\tau\ltimes(\mathrm{SO}(3) \times \mathrm{SO}(3))$. The left and right copies of the rotation group correspond respectively to rotations with respect to an inertial frame fixed in space and rotations with respect to an internal reference frame. The discrete factor corresponds to the transposition symmetry, which is a very important and famous feature of Dirichlet's problem. $G$ acts on itself as follows: If $(\gamma;g,h),(\gamma';g',h')\in \mathbb{Z}^\tau_2\ltimes(\mathrm{SO}(3)\times\mathrm{SO}(3))$ then
\begin{equation*} (\gamma;g,h)\cdot
(\gamma';g',h')=(\gamma\gamma';(g,h)\cdot(\gamma\cdot
(g',h')),\end{equation*}
 where for the nontrivial element $\sigma\in\mathbb{Z}^\tau_2$, $\sigma\cdot (g',h')=(h',g')$, and
$\mathrm{SO}(3)\times\mathrm{SO}(3)$ acts on itself by the direct product of left matrix multiplication. The  Lie algebra of $G$ is then $\g=\mathbb{R}^3\oplus\mathbb{R}^3$ on which $G$ acts by the
 adjoint representation given by
\begin{equation}\label{adjoint action}
\Ad_{(\gamma ; g,h)}
(\xi_L,\xi_R)=\gamma\cdot(g\cdot\xi_L,h\cdot\xi_R),
\end{equation}
where $\sigma\cdot (\xi_L,\xi_R)=(\xi_R,\xi_L)$ and $g\cdot\xi_L$ is
the rotation of $\xi_L$ by $g$ (and similarly for $\xi_R$). Analogously, we have for the coadjoint representation, after identifying $\g\simeq\g^*$ with the dot product,
\begin{equation}\label{coadjoint action}
 \Ad^*_{(\gamma ; g,h)^{-1}}
 (\mu_L,\mu_R)=\gamma\cdot(g\cdot\mu_L,h\cdot\mu_R).
\end{equation}

The action of $G$ on $\mathrm{SL}(3)$ is given by
\begin{equation}(e;L,R)\cdot F  =  LFR^T,\quad
(\sigma;L,R)\cdot F  =  RF^TL^T,\label{GactionSL3}
\end{equation}

Note that the $\mathbb{Z}^\tau_2$ transposition  symmetry on
$\mathrm{SL}(3)$ first noticed by Dedekind maps a rigidly rotating
configuration without internal motion  into one that is stationary but has internal vorticity.
These adjoint to each other ellipsoids are called Jacobi and Dedekind ellipsoids.
More generally, the transposition symmetry interchanges external
and internal rotations for any solution of Dirichlet's
model.

Dirichlet's problem can be formulated as a symmetric Hamiltonian system on the phase space $\PP = T^*\mathrm{SL}(3)$ with Hamiltonian
$$h(P_F)= \frac{1}{2T}\mathrm{tr}(P_F^TP_F)+V(F), \quad  P_F\in T^*\mathrm{SL}(3)$$ where the symplectic form $\omega$ is the canonical one on the cotangent bundle and the $G$-invariance of $\ll,\cdot,\cdot\gg$ and $V$ makes the natural lifted action of $G$ on $\PP$ a Hamiltonian one, with $G$-equivariant momentum map $\J:\PP\rightarrow \g^*$. Riemann ellipsoids (including those with axisymmetric configurations) correspond in this setting to relative equilibria of this Hamiltonian system. Recall that a relative equilibrium  of a symmetric Hamiltonian system is given by 
$$z(t) = e^{\xi t}\cdot z,$$
for some (possibly not unique) constant element $\xi\in \g_\mu = \mathrm{Lie}(G_\mu)$ called a velocity of the relative equilibrium. Here $G_\mu$ denotes the isotropy group, or stabilizer, of  $\mu=\J(z)\in\g^*$ with respect to the coadjoint representation of $G$. 

I will adopt the usual point of view of classifying relative equilibria by their stabilizers $G_z$. For Hamiltonian systems of mechanical type with phase space $T^*M$ like Dirichlet's problem, it is straightforward (\cite{ROSD02}) to characterize $G_z$ as

\begin{equation}G_{z}=\{g\in G_x\,:\,\Ad_g\xi-\xi\in\g_x\}.\label{stabREformula}\end{equation}

where $x=\pi_M(z)$, the projection of the phase space point to $M$ under the canonical cotangent bundle projection. In our problem $z=P_F$ and $x=F$. Using the singular value decomposition we have that every $F\in \mathrm{SL}(3)$ can be written as  $F=LAR^T$, with $A$ a diagonal matrix. This is equivalent to every element of the configuration space being in the $G$-orbit of a diagonal matrix under the action \eqref{GactionSL3}. As mentioned before, we consider the whole $G$-orbit as a single relative equilibrium and therefore we can study its properties taking a representative with diagonal configuration. Depending on the number of different eigenvalues of $A$ we can have a spherical (1), spheroidal (2) or ellipsoidal (3) configuration. Riemann ellipsoids are relative equilibria of the form 
\begin{equation}\label{defPF} P_F=\ll (\xi_L,\xi_R)_{\mathrm{SL}(3)}(F), \cdot\gg\end{equation} 
where $\xi_M$ denotes the fundamental vector field corresponding to the element $\xi\in\g$ under the $G$-action on $M$. In particular a MacLaurin spheroid is a solution in the $G$-orbit of $P_F$ corresponding to the data

\begin{equation}\label{diagF}F=\mathrm{diag}(a,a,c),\quad \xi_L,\xi_R\parallel \mathbf{e_3}, \quad \xi_L\neq \xi_R, \quad c<a.\end{equation}

The eccentricity of the spheroid is $e=\sqrt{1-c^2/a^2}$. The case $\xi_L=\xi_R$ would correspond to a spheroidal equilibrium, a solution that although is in principle possible for a Hamiltonian system with these configuration space and symmetries, is not present for the gravitational potential \eqref{gravpot} of Dirichlet's problem under the incompressibility condition. The precise parametrization of the MacLaurin family is obtained through the functional dependence of $\xi_L-\xi_R=\frac{\Omega(e)}{2}\mathbf{e_3}$ on $e$, given by the famous condition found by MacLaurin: 
\begin{equation}\label{MLformula}\frac{\Omega^2}{\pi \rho_0 G} = 2\frac{\sqrt{1-e^2}}{e^3} (3-2e^2)\arcsin(e) - \frac{6}{e^2} (1-e^2), \quad e\in(0,1).\end{equation}

For the case of the MacLaurin spheroid, using \eqref{stabREformula} is easy to obtain $G_z=\widetilde{\mathrm{O}(2)_\mathbf{e_3}}$, which is the group generated by the elements
$$(e;R_\mathbf{e_3}(\theta), R_\mathbf{e_3}(\theta)), \quad \mathrm{and}\quad (\sigma;R_\mathbf{n}(\pi),R_\mathbf{n}(\pi))$$
where $R_\mathbf{k}(\theta)$ denotes a rotation of angle $\theta\in [0,2\pi)$ around the vector $\mathbf{k}$ and $\mathbf{n}$ is a vector perpendicular to $\mathbf{e_3}$. Notice that for each value of the eccentricity many combinations of $\xi_L$ and $\xi_R$ producing the same solution are possible. This is a direct consequence of the fact that $G_z$ is continuous. The precise pattern of bifurcations, and particularly the fact that this family exhibits bifurcations without changes in its stability, is also a direct consequence of the continuity of $G_z$. 

In the following sections I will introduce the necessary material to apply the bifurcation analysis to the MacLaurin family. I will show how there are three possible different bifurcating branches of Riemann ellipsoids and how the analysis of the stabilizers of the bifurcating branches allows to the precise identification of these branches among the different types of Riemann ellipsoids.

\section{Hamiltonian bifurcations from parametrized branches of relative equilibria}\label{sec:Hambif}

Relative equilibria typically organize into clearly distinct families or branches which are parametrized by quantities like energy or angular momentum. 
In certain cases it happens that at specific values of the parameters two or more distinct families of relative equilibria intersect. In this situation we say that there is a bifurcation. These bifurcations could happen for isolated values of the parameters or for a continuous set of parameter values. The existence of such bifurcation points and bifurcating branches is intimately related to the changes in stability of the relative equilibria as the parameters are varied as well as to the various stabilizer groups present in the problem.
It is a standard fact that relative equilibria $z$ with velocity $\xi$ are characterized as  critical points via the equation 

\begin{equation}\label{firstvariation}\ed h_\xi(z)=0\end{equation}
where the augmented Hamiltonian function is defined as

 \begin{equation}\label{augHam} h_\xi (z)= h(z) - \langle \J(z), \xi \rangle \end{equation}

The problem of Hamiltonian bifurcations consists in, given a parametrized family $\bar z(w)$ of solutions of \eqref{firstvariation}, how the existence of other solutions near it depends on the parameters $w$.  This naturally leads to approach this problem within the field of critical point theory and more specifically to the study of the null spaces of the second variation of $h_\xi$ at its critical points in what sometimes has been referred to as ``energetic methods''. In this article I will adopt the convention that the bifurcating relative equilibria must be of a ``distinct type'' than the original equilibrium or branch. Otherwise we would refer to that situation as persistence of relative equilibria. We will use the following definition:
 We say that $\bar z: W\rightarrow \PP$ is a parametrized branch of relative equilibria of the same symplectic type if
 $W$ is an open set in a
vector space and $\bar z(w)$ is a relative equilibrium for all $w$ such that

\begin{itemize}
    \item[(i)] if $w_1\neq w_2$ then $\bar z(w_1)$ and $\bar z(w_2)$ belong to different group orbits, and
    \item[(ii)] for all $w$, all the groups $G_{\J(\bar z(w))}$ are conjugate and the same applies to the groups $G_{\bar z(w)}$.
\end{itemize}

Condition $(i)$ is necessary in order to reconcile this framework with reduction theory. Historically, a relative equilibrium has  been defined not as a point $z$ in phase space $\PP$ but rather as an equivalence class $[z]\in \PP/G$. Therefore, all different trajectories of the form $z(t)=e^{t\xi}\cdot z$ that belong  to the same $G$-orbit are considered to be the same relative equilibrium. 

Condition $(ii)$ is necessary for technical reasons. Intuitively, this represents the idea of all points in the branch being the ``same type'' or having the ``same shape''. The precise implementation of this notion of equivalence is through the comparison of the different symmetries available for each relative equilibrium. These are the stabilizers of the point $z$ itself and of its momentum value $\J(z)$. This second condition allows us to be precise about what is and what is not a member of a family of relative equilibria. For instance, it has been shown that among Riemann ellipsoids (see \cite{RoSD99} and \cite{ROSD07}), the spherical equilibrium satisfies $G_z = \mathbb{Z}^\tau_2 \times \mathrm{SO}(3)$ while as we have seen MacLaurin spheroids have stabilizers conjugate to $\widetilde{\mathrm{O}(2)_\mathbf{e_3}}$. For this reason, and in contrast with the approach taken in several authoritative works (see for instance \cite{Cha87}), I will not consider the spherical equilibrium as  ``the first member'' of the MacLaurin family. 

All the analysis in this article will follow from a slight modification to Theorem 5.2 in \cite{RO20}. Before stating it, recall that if the $G$-action on $\PP$ is proper, for any $z\in\PP$ with momentum $\J(z)=\mu$ the stabilizer $G_z$ must be compact and, because of the equivariance of $\J$, $G_z\subseteq G_\mu$. Therefore there exists a $G_z$-invariant splitting
\begin{equation}\label{gmusplitting}\g_\mu = \m \oplus \g_z   \end{equation}
We then have
\begin{thm}\label{thm bifurcations} Let $z\in\PP$ be a relative equilibrium with momentum $\J(z)=\mu$ satisfying $[\g_\mu, \g_\mu]=0$ and velocity
$\xi\in\g_\mu$ written as $\xi=\xi^\perp+\eta$ according to \eqref{gmusplitting}.
Let $W\subseteq {\m^*}^{G_z}
\times \g_z^{(G_z)_\eta}$ be a linear subspace and $\bar z:W\rightarrow \PP$
 a parametrized branch of relative equilibria of the same symplectic type,
satisfying $\bar z(0,0)=z$, $\J(\bar z(\rho,\eta'))=\mu+\rho$,
$G_{\J(\bar z(\rho,\eta'))}=G_\mu$ and $G_{\bar z(\rho, \eta')}=G_z$. Let
$\xi(\rho, \eta')$ be a family of velocities for points of this
branch chosen such that $\Proj_{\g_z}(\xi(\rho, \eta'))=\eta + \eta'$. Now suppose that for some $L\subseteq (G_z)_\eta$
\begin{itemize}
    \item[(i)] The normalizer $N_{(G_z)_\eta}(L)$ acts on $\Null:=\ker \ed^2_zh_{\xi}\rrestr{N^L}$ with cohomogeneity one,   and
    \item[(ii)] the eigenvalue $\sigma(\rho, \eta')$ of $\ed^2_{\bar z(\rho, \eta')}h_{\xi(\rho, \eta')}\rrestr{\Null}$ crosses $0$ at $0\in W$.
\end{itemize}
Then, for every  $v\in \Null$ close enough to the origin,  there is a relative equilibrium
$\zt_v$ near $z$ and
not in the $G$-orbit of any point of the original branch $\bar z(\rho, \eta')$.
This point $\zt_v$ satisfies $G_{\zt_v} =  (G_z)_v$.
\end{thm}

\noindent{{\bf Remarks.}}
\begin{itemize}
\item[(i)] In \cite{RO20} this result is stated for Abelian momentum isotropy groups, but the proof and results translate with minor changes to the weaker case where only their Lie algebras are Abelian (Section 2.3 and Theorem 5.1 in \cite{MoRO16}).

\item[(ii)] The cohomogeneity condition in (i) is always satisfied in the cases that will appear in the study of the MacLaurin spheroid since the spaces $\Null$ will be  either one-dimensional or two-dimensional supporting a representation of $SO(2)$ with a one dimensional orbit space, which makes the problem equivalent to a one-dimensional one. 

\item[(iii)] $N$ refers to the symplectic normal space at $z$, defined as the quotient space $\ker T_z\J/\g_\mu\cdot z$. It can be realized as any $G_z$-invariant complement to $\g_\mu\cdot z = T_z(G_\mu\cdot z)$ in $\ker T_z\J\subset T_z\PP$. The restriction $\ed^2_zh_{\xi}\vert_N$ is usually called the stability form due to the role it plays in the study of stability of relative equilibria using the Energy-Momentum method (see \cite{Pat92}).
\end{itemize}

\section{Bifurcations from the MacLaurin family}\label{secbif}
I will now apply the framework introduced in the previous sections to the study of Riemann ellipsoids near the MacLaurin family.
Most of the technical ingredients for the application to the  MacLaurin spheroid of the methods of Section \ref{sec:Hambif} and in particular of Theorem \ref{thm bifurcations} have been computed in \cite{ROSD07}. In order to avoid lengthy and redundant derivations I will rely on many of the computations in that reference. The first two paragraphs of this section follow closely the terminology of \cite{RoSD99}.
\ \\
 
\paragraph{\bf Riemann's Theorem.}
Riemann's work on Dirichlet's problem \cite{Rie61} provides necessary conditions on its steady solutions (relative equilibria in the symmetric Hamiltonian context) stated in the famous Riemann's theorem: A solution for Dirichlet's problem in which the length of the semiaxes as well as the angular velocity and vorticity vectors remain constant must be of one of the following types:

\begin{itemize}
\item[(i)] For a spherical configuration, $\xi_L$ and $\xi_R$ must be parallel.
\item[(ii)] For a spheroidal configuration $\xi_L$ and $\xi_R$ must lie in a plane containing the symmetry axis.
\item[(iii)] For a ellipsoidal configuration either $\xi_L$ and $\xi_R$ are parallel to the same principal axis of the ellipsoid or they lie in the same principal plane.
\end{itemize}

Historically, solutions for which $\xi_L$ and $\xi_R$ are parallel to a principal axis are called $S$-ellipsoids. Otherwise they are classified in  types $I$, $II$ or $III$. The distinction of an ellipsoid with three different semiaxes $(a_1, a_2, a_3)$ into the three different types is as follows (assuming $\xi_L$ and $\xi_R$ belong to the plane $ij$):

\begin{itemize}
 \item type $I$ if $a_k\geq \frac 12 (a_i + a_j)$.
 \item type $II$ if $a_k\leq \frac 12 \vert a_i - a_j\vert$ where $a_k$ is the middle axis.
 \item type $III$ if $a_k\leq \frac 12 \vert a_i - a_j\vert$ where $a_k$ is the smallest axis. 
\end{itemize}

\noindent{\bf Remark.}
The division of solutions with ellipsoidal configuration into ellipsoids of type I, II and III is not directly related to the symmetries of the problem. However it is important to notice that in this article we will not find any ellipsoid of types II or III, since it is clear from their definitions that the only ellipsoidal shapes that can occur near a spheroid are necessarily of type I.
\ \\
\paragraph{\bf Symmetries of Riemann ellipsoids and conditions on the possible bifurcations.}
As noted in the previous section, one possibility to classify families of relative equilibria is by the conjugacy class of their stabilizers $G_z$. Basic considerations on the topology of the $G$-action on $\mathrm{SL}(3)$ imply that possible Riemann ellipsoids bifurcating from the MacLaurin family need to have stabilizers conjugated to proper subgroups of the stabilizer of the MacLaurin family $\widetilde{\mathrm{O}(2)}$. Moreover, the stabilizer $G_z=G_{P_F}$ must be contained in $G_x=G_F$ as can be seen from \eqref{stabREformula}. Recall that there are 3 possible distinct configurations up to conjugacy, and from \eqref{GactionSL3} we get for each one:
\begin{itemize}
\item[(i)] Spherical configuration: $F=\mathrm{diag}(1, 1, 1), \quad G_F=\mathbb{Z}_2^\tau\times \mathrm{SO}(3)^D$
\item[(ii)] Spheroidal configuration: $F=\mathrm{diag}(a, a, c), \quad G_F=\mathbb{Z}_2^\tau\times \mathrm{O}(2)_\mathbf{e_3}^D$
\item[(i)] Ellipsoidal configuration: $F=\mathrm{diag}(a, b, c), \quad G_F=\mathbb{Z}_2^\tau\times {\mathbb{D}_2}_\mathbf{e_3}^D$
\end{itemize}
with $a\neq b\neq c$ and $\mathrm{det}(F)=1$. Here $K^D$ denotes the diagonal embedding of $K\subset G$ into $G\times G$. $\mathrm{O}(2)_\mathbf{e_3}^D$ is generated by $R_\mathbf{e_3}(\theta)$ and $R_\mathbf{n}(\pi)$ while ${\mathbb{D}_2}_\mathbf{e_3}^D$ is generated by $R_\mathbf{n}(\theta)$ and $R_\mathbf{n}(\pi)$. From the above, we find that the topological restrictions on the possible bifurcations at the configuration level imply that the only relative equilibria that can bifurcate from a branch with spheroidal configurations must have spheroidal or ellipsoidal configurations. The existence conditions and stability properties for all possible spherical or spheroidal configurations of Dirichlet's problem have been studied in \cite{ROSD07} 
using the same geometric formalism as this article. There it is shown that the only possible Riemann ellipsoids with spheroidal configurations are the MacLaurin spheroids and the transversal (or prolate) spheroids, which are contained in the family of S ellipsoids (see also 48.(e) in \cite{Cha87}). It follows from any of these two references that the transversal family cannot be joined with the MacLaurin family and therefore our restrictions imply that the only possible bifurcating solutions from the MacLaurin family must have exclusively ellipsoidal configurations. This fact will allow to identify the bifurcating Riemann ellipsoids obtained by Theorem \ref{thm bifurcations} by comparing their predicted stabilizers with the classification of all possible relative equilibria and their stabilizers compatible with Dirichlet's problem carried out in \cite{RoSD99}.
\ \\
\paragraph{\bf Characterization of the MacLaurin branch.} In the remainder of this article we will use the following setup: we can rewrite \eqref{diagF} in terms of the eccentricity and obtain the configuration of a MacLaurin spheroid of  as
\begin{equation}\label{confML}
F=\mathrm{diag}((1-e^2)^\frac{-1}{6}, (1-e^2)^\frac{-1}{6},  (1-e^2)^\frac{1}{3})
\end{equation}
From the splitting \eqref{gmusplitting} we have $\g_z=\langle(\mathbf{e_3},\mathbf{e_3})\rangle$ and $\m^*=\langle(\mathbf{e_3},-\mathbf{e_3})\rangle$ once we identify $\m$ and $\m^*$ via the standard inner product in $\mathbb{R}^3$. 
A possible way of choosing the velocity $(\xi_L,\xi_R)$ for a MacLaurin spheroid of eccentricity $e\in (0,1)$  according to Theorem \ref{thm bifurcations} is 
\begin{equation}\xi^\perp  =  \frac{\Omega(e)}{2}(\e_3,-\e_3), \quad \eta = 0 \label{velML}
    \end{equation}
Here $\Omega (e)$ is given by equation \eqref{MLformula}. Along the MacLaurin family we have the momentum map, expressed as the angular momentum/circulation pair given by
    \begin{equation}\J(z)=(\mathbf{j},\mathbf{c})=2(1-e^2)^{-1/3}T\Omega(\mathbf{e_3},-\mathbf{e_3}) = \hat\mu(e)  (\mathbf{e_3},-\mathbf{e_3})=\mu \label{muML}\end{equation}

\begin{figure}
    \centering
    \includegraphics[width=0.5\linewidth]{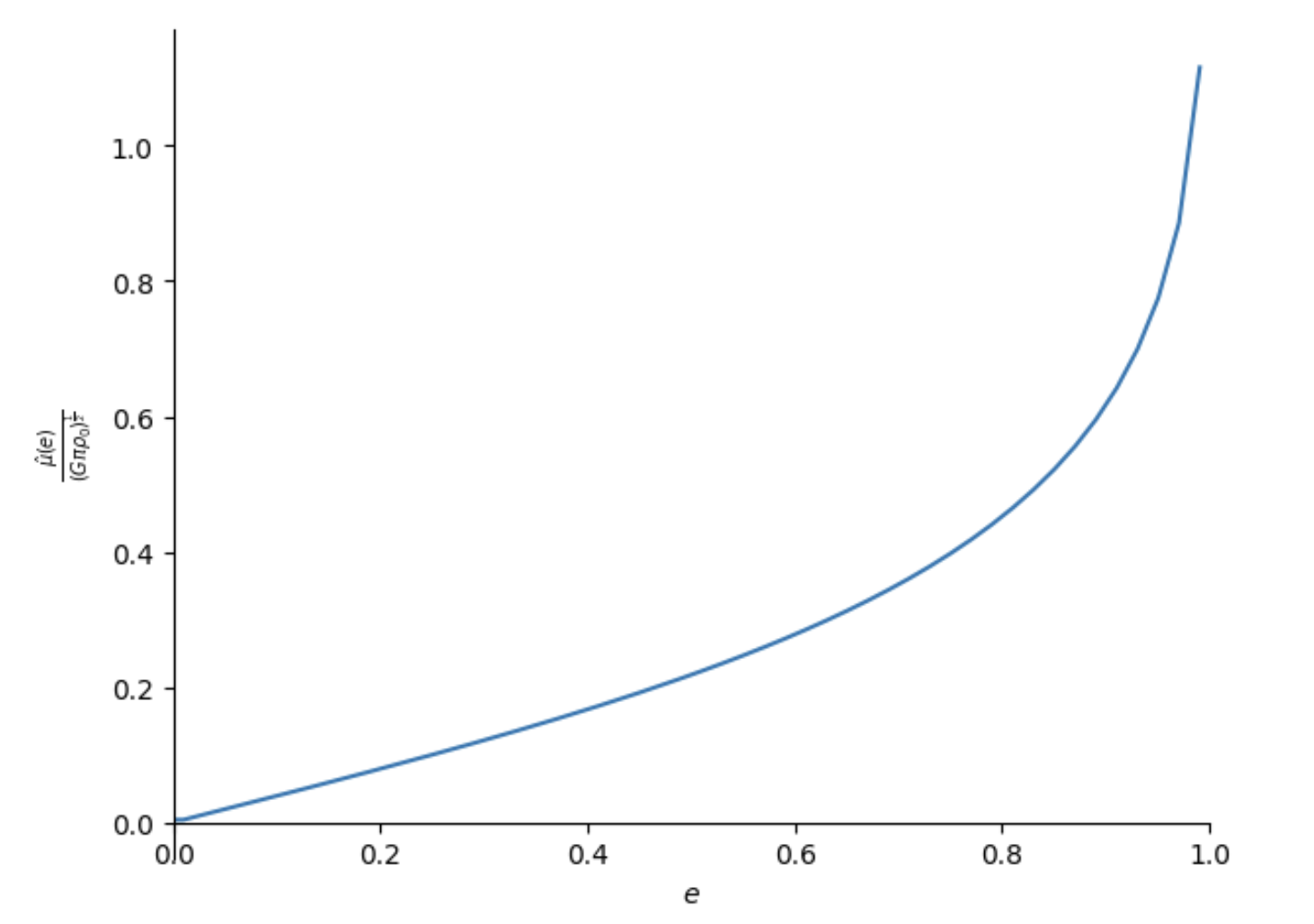}
    \caption{Graph of $\hat \mu(e)$ in $(G\pi\rho_0)^\frac{1}{2}$ units along the MacLaurin family. It allows parametrizing the MacLaurin family by $e\in(0,1)$ or by $\hat\mu\in (0,\infty)$.  }
    \label{fig:mu(e)}
\end{figure}

It follows that  $G_\mu=\widetilde{\mathrm{SO}(2)_{\mathbf{e_3}}\times\mathrm{SO}(2)_{\mathbf{e_3}}}$ which is generated
by the elements $(\sigma;R_\mathbf{n}(\pi),R_\mathbf{n}(\pi))$ and $(e;R_\mathbf{e_3}(\theta_1),R_\mathbf{e_3}(\theta_2))$ with $\mathbf{n}\perp \mathbf{e_3}$ and $\theta_1,\theta_2$ arbitrary. Along the family of MacLaurin spheroids the function $\hat \mu(e)=2(1-e^2)^{-1/3}T\Omega$ provides an injective map from the set of conjugacy classes of MacLaurin spheroids into $\m^*$ given by 
    \begin{equation}\label{invertMLbranch}e\mapsto \hat \mu(e)(\mathbf{e_3},-\mathbf{e_3})
    \end{equation} 
which is just the momentum value of each spheroid (see Figure \ref{fig:mu(e)}). Later we will use this map in order to define the branch $\bar z(\rho,\eta')$ of Theorem \ref{thm bifurcations}.
\ \\
\paragraph{\bf The symplectic normal space}
The symplectic normal space at a MacLaurin spheroid can be expressed as $N\simeq \q^\mu\oplus\mathbf{S}\oplus \mathbf{S}^*$ on which $G_z$ acts diagonally and where
$$
\q^\mu  =  \langle t_1,t_2,t_3,t_4\rangle$$
with 
\begin{equation}\label{qumubasis} t_1=(\mathbf{e_1},0),\quad t_2=(0,\mathbf{e_1}),\quad t_3=(\mathbf{e_2},0),\quad t_4=(0,\mathbf{e_2})\end{equation}
and $\mathbf{S}$ can be identified with elements of the form $(a_1,a_2,a_3)\in\R^3$ via the map 

\begin{equation}\label{Sbasis}
    (a_1, a_2, a_3) \mapsto \left( \begin{array}{ccc} 
        a_1+a_2 & a_3       & 0 \\          
        a_3     & a_1-a_2   & 0 \\
        0       & 0         & -2\frac ca a_1 \end{array}\right)
\end{equation}
(see equation (51) in \cite{ROSD07}). 
According to the previous splitting we have the following expression for the symplectic matrix  of $N$.
    $$\Omega_N=\left(\begin{array}{ccc} \Xi & 0 & 0 \\ 0 & 0 & I \\ 0 & -I & 0\end{array}\right)$$
with $I$ the $3\times 3$ identity matrix and
    $$\Xi=\frac{2T\Omega}{(1-e^2)^{1/3}}\left(
        \begin{array}{cccc} 
            0 & 0 & -1 & 0 \\ 
            0 & 0 & 0 & 1  \\ 
            1 & 0 & 0 & 0  \\ 
            0 & -1 & 0 & 0 
        \end{array}\right)$$
\ \\
\paragraph{\bf The action and momentum map.}
We can choose two generators for  $G_z=\widetilde{\mathrm{O}(2)_\mathbf{e_3}}$ as 
    \begin{eqnarray*} 
        \theta     & = & (e;R_\mathbf{e_3}(\theta)R_\mathbf{e_3}(\theta))\\
        \bar\sigma & = & (\sigma;R_\mathbf{e_1}(\pi), R_\mathbf{e_1}(\pi))
    \end{eqnarray*}
where $R_\mathbf{e_1}(\pi)=\mathrm{diag}(1,-1,-1)$. On each factor of $N$ we have: 
\begin{itemize}
    \item[(i)] On $\mathbf{S}$ and $\mathbf{S}^*$: From \eqref{Sbasis} it is clear that the linearization of the action \eqref{GactionSL3} produces 
        \begin{eqnarray} 
            \theta\rrestr{\mathbf{S}}       & = & \left(\begin{array}{ccc} 
            1   & 0                 & 0 \\             
            0   & \cos (2\theta)    & -\sin (2\theta) \\ 
            0   & \sin(2\theta)    & \cos (2\theta) \end{array}\right)   \label{thetaS}\\
            \bar\sigma\rrestr{\mathbf{S}}   & = & \mathrm{diag}(1,1,-1) \label{sigmaS}
        \end{eqnarray}
and the same expressions work for $\Sl^*$ under the isomorphism induced by the trace inner product \eqref{metricsl3}.
    \item[(ii)] On $\q^\mu$: Similarly, from \eqref{qumubasis} it is straightforward to obtain the restriction to $\q^ \mu$ of the adjoint action \eqref{adjoint action}
        \begin{eqnarray} 
        \theta\rrestr{\q^\mu}       & = & \left(\begin{array}{cccc} 
                                                    \cos (\theta)   & 0             & -\sin(\theta)     & 0 \\ 
                                                             0      & \cos (\theta) & 0                 & -\sin(\theta) \\                   
                                                    \sin(\theta)    & 0             & \cos(\theta)      & 0   \\ 
                                                          0         &  \sin(\theta) &  0                & \cos(\theta) \end{array}\right)     \label{thetaqmu}\\ 
        \bar\sigma\rrestr{\q^\mu}   & = &  \left(\begin{array}{cccc} 
                                                0 & 1 & 0 & 0 \\ 
                                                1 & 0 & 0 & 0 \\ 
                                                0 & 0 & 0 & -1 \\ 
                                                0 & 0 & -1 & 0\end{array}\right)\label{sigmaqmu}    
        \end{eqnarray}

    \item[(iii)] On $\m^*$: Since $\m^*=\langle (\mathbf{e_3},-\mathbf{e_3})\rangle$ it follows from \eqref{coadjoint action} that the action of $G_z$ is trivial on $\m^*$ and then ${\m^*}^{G_z}=\m^*$.
\end{itemize}
Let $\eta\in\R$ and identify it with the element $\eta(\mathbf{e_3},\mathbf{e_3})\in\g_z$. The corresponding infinitesimal generators are

    $$\eta\rrestr{\Sl}  =  2\eta\left(\begin{array}{ccc} 
        0 & 0 & 0 \\ 
        0 & 0 & -1 \\ 
        0 & 1 & 0   \end{array}\right)
    \quad \eta\rrestr{\q^\mu}  =  \eta\left(\begin{array}{cccc} 
        0 & 0  & -1 & 0 \\ 
        0 & 0  & 0  & -1 \\ 
        1 & 0  & 0  & 0 \\ 
        0 & 1  & 0  & 0        \end{array}\right) $$
With these expressions we can compute the quadratic momentum map for the $G_z$ action on  $(N,\Omega_N)$, with the standard form 
$$\mathbf{J}^\eta_N(v) = \frac 12 \Omega_N(\eta\cdot v, v)$$ 
where $\eta\in\g_z$, $v\in N$ and $\eta\cdot v$ denotes the infinitesimal action on $N$. It is straightforward to obtain in our case 
\begin{equation}\label{momentumslice}\J_N^\eta(\lambda;a,\beta) = \frac{T\Omega\eta}{(1-e^2)^{1/3}}(\lambda_1^2-\lambda_2^2+\lambda_3^2-\lambda_4^2)+2\eta(a_2\beta_3-a_3\beta_2)\end{equation}
\ \\
\paragraph{\bf The stability form.}
The Hessian of the augmented Hamiltonian for $\xi = \xi^\perp$ is obtained in Section 6.2 of \cite{ROSD07} and block-diagonalizes on the product $\q^\mu\times \Sl \times \Sl^*$ as

    $$\mathbf{d}^2_zh_{\xi^\perp}\rrestr{N}=\left(\begin{array}{ccc} Ar & 0   & 0 \\ 
                                                                0 & R_2 & 0 \\ 
                                                                0 & 0   & R^{-1}_1
                                            \end{array}\right)$$
where the block $\mathrm{Ar}$ is called the Arnold form and has the expression 

    \begin{eqnarray*} 
        Ar & = & \left(\begin{array}{cccc} 
                                           A_1 & -A_2 & 0    & 0 \\ 
                                          -A_2 & A_1  & 0    & 0 \\ 
                                             0 & 0    & A_1  & -A_2 \\ 
                                             0 & 0    & -A_2 & A_1
                        \end{array}
                  \right)
    \end{eqnarray*}    
with

    \begin{eqnarray*}
        A_1 & = & \frac{(8-e^4 -4e^2)T\Omega^2}{e^4(1-e^2)^{\frac 13}}\\ 
        A_2 & = & \frac{8(1-e^2)^{\frac 16} T\Omega^2}{e^4}
    \end{eqnarray*} 
The blocks on $\Sl$ and $\Sl^*$ are diagonal and of the form

    \begin{eqnarray*}
        R_1 & = & \diag(2T(3-2e^2),2T,2T)\\ 
        R_2 & = & \diag(S_1,S_2,S_2)
    \end{eqnarray*}
where

    \begin{eqnarray}
        S_1 & = & \frac{2R}{e^5}\left(9e(3-5e^2+2e^4)-\sqrt{1-e^2}(27-36e^2+8e^4)\arcsin(e)\right)\\
        S_2 & = & \frac{R}{e^5}\left(e (1-e^2)(3+4e^2)-\sqrt{1-e^2}(3+2e^2-4e^4)\arcsin(e)\right)\label{S2}
    \end{eqnarray} 

From \eqref{momentumslice} it is immediate to compute the second variation of $\J_N^\eta$ obtaining

    $$\mathbf{D}_z^2\J_N^\eta = \left(\begin{array}{ccc} 
            \eta(1-e^2)^{-\frac 13}T\Omega J_1     & 0          & 0          \\ 
            0                           &            & 2\eta J_2 \\ 
            0                           & -2\eta J_2 & 0
    \end{array}\right)$$
with $J_1=\mathrm{diag}(1,-1,1,-1)$ and  
    $$J_2=\left(\begin{array}{ccc} 
            0 & 0 & 0  \\ 
            0 & 0 & 1 \\ 
            0 & -1 & 0
                       \end{array}
          \right)$$

Let $z$ be a phase space point that corresponds to a MacLaurin spheroid with eccentricity $e$ and velocity $\xi=\xi^\perp+\eta$ with $\xi$ given in \eqref{velML} and $\eta\in\g_z$. Notice that by construction $$\ed^2_{z}h_{\xi^\perp+\eta}\vert_N=0$$
for any value of $\eta$. It follows from the previous discussion and \eqref{augHam} that
\begin{equation}\label{augHessN}\ed^2_{z}h_{\xi^\perp+\eta}\vert_N 
    =\ed^2_{z}h_{\xi^\perp}\vert_N 
    - \mathbf{D}_z^2\J_N^\eta
    =\left(\begin{array}{ccc} 
            Ar-\eta(1-e^2)^{-\frac 13}T\Omega J_1 & 0           & 0          \\ 
            0                         & R_2         & -2\eta J_2 \\ 
            0                         & 2\eta J_2    & R_1^{-1}
    \end{array}\right)
    \end{equation}
Notice that the matrix \eqref{augHessN} is block-diagonal. We will now study its eigenvalues on the blocks $\q^\mu$ and $\Sl \times \Sl^*$ which in turn will produce different kinds of bifurcations as we shall see in a later paragraph.\\

\begin{itemize}
\item[(i)] On $\mathbf{S}\times \mathbf{S^*}$: The matrix 
$$\left(\begin{array}{cc} R_2 & -2\eta J_2\\  2\eta J_2 & R_1^{-1} \end{array}\right)$$
has eigenvalues
    \begin{eqnarray*}
        \sigma^2_a      &= & S_1\\ 
        \sigma^2_b      & = &\frac{1}{2T(3-2e^2)}\\ 
        \sigma^2_{\pm}  & = & \frac 12\left(\frac{1}{2T}+S_2\pm\sqrt{\left(\frac{1}{2T}-S_2\right)^2+16\eta^2}\right)\end{eqnarray*}

\item[(ii)] On $\mathfrak{q}^\mu$: Notice that 
$$Ar-\eta(1-e^2)^{-\frac 13}T\Omega J_1=\left(\begin{array}{cc} B & 0 \\ 0 & B\end{array}\right)$$
where
$$B=\left(\begin{array}{cc} A_1-\eta(1-e^2)^{-\frac 13}T\Omega & -A_2 \\ -A_2 & A_1+\eta(1-e^2)^{-\frac 13}T\Omega \end{array}\right)$$
with eigenvalues
\begin{equation*}\sigma^1_\pm=A_1\pm \sqrt{A_2^2+\eta^2(1-e^2)^{-\frac 23}T^2\Omega^2}\end{equation*}
\end{itemize}
\ \\
\paragraph{\bf Zero eigenvalues.}
It is clear that $\sigma^1_+,\, \sigma^2_+,\, \sigma^2_a$ and $\sigma^2_b$ are always positive. Also, 
 for any value of $e$, $A_1>A_2>0$, so $\sigma^1_-$ is a zero eigenvalue when  
    $$\eta^2=\eta_1^2(e)=\frac{A_1^2-A_2^2}{T^2\Omega^2}(1-e^2)^\frac 23$$
The only other remaining possibility for a zero eigenvalue is when $\sigma^2_-=0$. We can distinguish two cases:
\begin{enumerate}
\item $S_2\geq 0$. Then $\sigma^2_-=0$ when
$$\eta^2=\eta^2_2(e)=\frac{S_2}{8T}.$$
\item $S_2< 0$. Then $\sigma^2_-<0$ for any value of $\eta$ since $T$ is a positive constant.
\end{enumerate}

\begin{figure}
    \centering
    \includegraphics[width=0.7\linewidth]{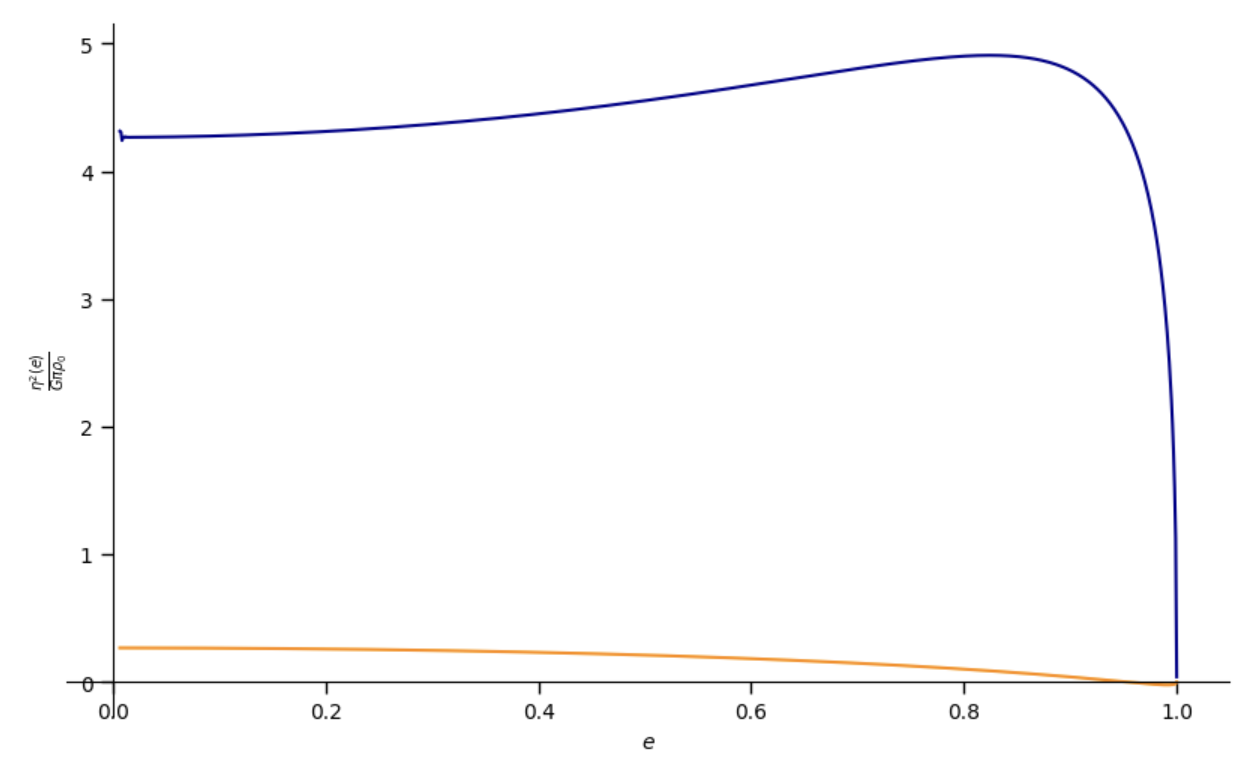}
    \caption{Graphs of $\eta_1^2(e)$ (top) and $\eta_2^2(e)$ (bottom) along the MacLaurin family expressed in $G\pi\rho_0$ units. Notice how $\eta_2$ becomes imaginary at $e_\mathrm{crit}\simeq 0.952887$, which is the point where the MacLaurin spheroid becomes unstable, while $\eta_1(e)$ is always real. There is no value of the eccentricity for which $\eta_1(e)=\eta_2(e)$.}
    \label{fig:eta2(e)}
\end{figure}

We will call $e_\mathrm{crit}$ the value for which $S_2=0$. By \eqref{S2} we have $e_\mathrm{crit}\simeq 0.952887$. If $e<e_\mathrm{crit}$ then $S>0$ and the MacLaurin spheroid is  stable, if $e>e_\mathrm{crit}$ then $S<0$ and the MacLaurin spheroid is unstable (see Figure \ref{fig:eta2(e)} and Theorem 6.2 in \cite{ROSD07}).

In order to apply Theorem \ref{thm bifurcations} let us consider the map  $\mu(e)=\hat\mu(e)  (\mathbf{e_3},-\mathbf{e_3})$ given in \eqref{muML}. We will be interested in the possible bifurcations that can occur at a given element of the family $z_0=z(e_0)$ with momentum $\mu_0=\mu(e_0)$.  We can use \eqref{invertMLbranch} to reparametrize $\rho\in \m^*$ by $e$ with the expression $\rho(e) = \mu(e) -\mu_0$ which can be inverted due to the invertibility of $\hat\mu(e)$. Therefore the function $\bar z(\rho) = z(e)$ satisfies $\bar z(0)= z_0$ and consists again on MacLaurin spheroids of the same symplectic type since
    \begin{eqnarray*}
    \J(\bar z(\rho)) = \mu(e(\rho)) & = & \mu_0 + \rho  \\
    G_{\J(\bar z(\rho))}      & = & G_{\mu_0} \\
    G_{\bar z(\rho)} & = & G_{z_0}
    \end{eqnarray*}
The second identity follows since $\rho\parallel\mu_0$. The last identity is a consequence of \eqref{stabREformula} and \eqref{defPF}
with $F$ and $\xi^\perp$  given by \eqref{confML} and \eqref{velML}. 

At $e_0$ we find that the kernels of $\ed^2_{z(e_0)}h_{\xi^\perp+\eta_1(e_0)}\vert_N$, and $\ed^2_{z(e_0)}h_{\xi^\perp+\eta_2(e_0)}\vert_N$ are respectively
\begin{eqnarray*}
N_0(\sigma_-^1) & = & \langle (D,1,0,0),(0,0,D,1) \rangle\subset \q^\mu\\
N_0(\sigma_-^2) & = & \langle (0,0,1,0,2T\eta_2(e_0),0),(0,1,0,0,0,-2T\eta_2(e_0))\rangle\subset \Sl\times \Sl^*
\end{eqnarray*}
where
\begin{equation}\label{D}D=\frac{A_1-\eta_1(e_0)(1-e_0^2)^{-\frac 13}T\Omega}{A_2} \end{equation} 
We can address three different cases:
\ \\
\paragraph{\bf Type I ellipsoids.}
Consider the family of MacLaurin spheroids given by $\bar z(\rho)$ and 
$$\xi(\rho, \eta')= \xi^\perp(e(\rho)) +\eta_1(e_0)+ \eta'$$ 
for $\eta'\in\g_z$. Recall that we are identifying $\eta\in \mathbb{R}$ with $\eta({\bf e_3}, {\bf e_3})\in\g_z$. Since $\eta_1(e_0)\neq 0$ we find, using \eqref{adjoint action}, that $(G_z)_{\eta_1(e_0)}=\mathrm{SO}(2)_\mathbf{e_3}^D\subset \mathrm{O}(2)_\mathbf{e_3}^D$. Let $L=\{e\}$ in Theorem \ref{thm bifurcations}. Then $W=\m^*\times \g_z$ 
and $N_{\mathrm{SO}(2)_\mathbf{e_3}^D}(L)=\mathrm{SO}(2)_\mathbf{e_3}^D$.
From \eqref{thetaqmu} we see that its action on $N_0(\sigma_-^1)^{\{e\}}=N_0(\sigma_-^1)$ is equivalent to the standard representation of $S^1$ on $\mathbb{R}^2$. Therefore, this action is of cohomogeneity 1. The eigenvalue $\sigma_-^1$ crosses 0 at $(0,0)\in W$.
Also, it is easy to see  from \eqref{D} that $D\neq 1$ for any $e_0\in (0,1)$ which together with \eqref{sigmaqmu} implies that $\bar\sigma$ cannot fix any element of $N_0(\sigma_-^1)^{\{e\}}$. From this property and the form of \eqref{thetaqmu} we have $(G_z)_v = \{e\} = \mathbbm{1}$ if $v\neq 0$ and it follows that the bifurcating branch $\tilde z:N_0(\sigma_-^1) \rightarrow T^*\mathrm{SL}(3)$ predicted by Theorem \ref{thm bifurcations} satisfies 
$G_{\tilde z_v}=\mathbbm{1}$ for every $v\neq 0$. By the classification of the stabilizers obtained in \cite{RoSD99} we have that this branch consists of type I ellipsoids given by ellipsoidal configurations
that bifurcate from the MacLaurin family at every value $e_0\in(0,1)$. 
\ \\
\paragraph{\bf Type $S$ ellipsoids.}
Consider the family of MacLaurin spheroids given by $\bar z(\rho)$ and 
$$\xi(\rho, \eta')= \xi^\perp(e(\rho)) +\eta_2(e_0)+ \eta'$$ 
We will assume that $e_0\in (0, e_\mathrm{crit})$ which implies that $\eta_2(e_0)\neq 0$ and again  $(G_z)_{\eta_2(e_0)}=\mathrm{SO}(2)_\mathbf{e_3}^D$. We choose $L=\Z_{2\mathbf{e_3}}^D\subset {\mathbb{D}_{2}}_\mathbf{e_3}^D$ and 
then $N_0(\sigma_-^2)^{\Z_{2\mathbf{e_3}}^D}=N_0(\sigma_-^2)$. We have again $N_{\mathrm{SO}(2)_\mathbf{e_3}^D}(L)=\mathrm{SO}(2)_\mathbf{e_3}^D$ and by \eqref{thetaS} its action on $N_0(\sigma_-^2)$ is equivalent to the representation of weight 2 given by

$$\theta \mapsto    \left( 
    \begin{array}{cc}
    \mathrm{cos}(2\theta) & -\mathrm{sin}(2\theta) \\
    \mathrm{sin}(2\theta) & \mathrm{cos}(2\theta)
    \end{array}
                    \right)
$$
from which we see that the cohomogeneity condition is satisfied. The parameter space is given by $W=\m^*\times \g_z$ and $\sigma_-^2$ crosses 0 at $(0,0)\in W$. It follows from \eqref{sigmaS} that $\bar\sigma$ cannot fix any point on $N_0(\sigma_-^2)$. This together with \eqref{thetaS} imply that the bifurcating branch predicted by Theorem \ref{thm bifurcations} satisfies $G_{\tilde z_v}=(G_z)_v = {\mathbb{Z}_2}^D_\mathbf{e_3}$ if $v\neq 0$ and therefore it  consists of type $S$ ellipsoids where $\xi_L,\xi_R\neq 0$ are parallel and of different lengths as obtained in \cite{RoSD99}.
\ \\
\paragraph{\bf Type adjoint $S$ ellipsoids.}\label{sec:adjointS}
This case happens only at $e_0=e_\mathrm{crit}$ for which $\eta_2(0)=0$. Therefore $(G_z)_{\eta_2(0)} = G_z =\widetilde{\mathrm{O}(2)_\mathbf{e_3}}$. We consider the subgroup $L=\widetilde{\mathbb{D}_{2\mathbf{e_3}}}\subset \widetilde{\mathrm{O}(2)_\mathbf{e_3}}$, the parameter space  $W=\m^*$ and  then
$$\mathcal{N}=N_0(\sigma_-^2)^{\widetilde{\mathbb{D}_{2\mathbf{e_3}}}}=
\langle (0,1,0,0,0,0)\rangle\subset \Sl\times \Sl^*$$
Since $\dim\,\mathcal{N}=1$ the cohomogeneity one condition is satisfied and, just as in the two previous cases, $\sigma_-^2$ crosses $0$ at $0\in W$. Similar considerations as in the previous two cases show that for every $0\neq v\in\mathcal{N}$ we have $(G_z)_{\tilde z_v}=(G_z)_v=\widetilde{\mathbb{D}_{2\mathbf{e_3}}}$ and the bifurcating family consists on adjoint $S$ ellipsoids according to \cite{RoSD99}. It is easy to see that these are a particular kind of S ellipsoids for which $\xi_L=-\xi_R\neq 0$.

\noindent{\bf Remarks.}

\begin{itemize}
 
    \item[(i)] There are no other bifurcations from the MacLaurin family than the ones obtained here. The reason is that there is no other velocity field $\xi(\rho, \eta')$ for which 
    $$\ed^2_{\bar z(\rho)}h_{\xi(\rho, \eta')}\vert_{N}$$
    is degenerate at $\rho = 0$ and a straightforward application of the implicit function theorem implies that no more bifurcations can be produced. 

    \item[(ii)] The result that type $I$ ellipsoids bifurcate at any point in the MacLaurin family but bifurcations to $S$ ellipsoids only occur in the stable range of eccentricities $(0,e_\mathrm{crit})$ is well known and already appears in \cite{Cha87} although its linear analysis produces necessary but nor sufficient conditions for these bifurcations to occur.

    \item[(iii)] The adjoint $S$ ellipsoids found here form the second of the two bounding families of adjoint $S$ solutions reported by Chandrasekhar in Section 48.(c) of \cite{Cha87}. The other one bifurcates from the spherical equilibrium, regarded by the author as the first member of the MacLaurin family, a point of view that I don't adopt here based on symmetry reasons as explained in Section \ref{sec:Hambif}. These two families appear then at values $e=0,e_\mathrm{crit}$ and for all values in between the MacLaurin spheroid is stable and bifurcates to type $S$ ellipsoids.

    \item[(iv)] The Jacobi and Dedekind ellipsoids are famous solutions of type $S$ consisting on a rigidly rotating ellipsoid without internal motions or a static ellipsoid where the only rotations are internal respectively. Of course in the setup of modern Hamiltonian dynamics these are just two equivalent solutions since they belong to the same group orbit. The pair $(\xi_L,\xi_R)$ must be $(\xi,0)$ for the Jacobi ellipsoid and $(0,\xi)$ for the Dedekind ellipsoid with $\xi\parallel\mathbf{e_3}$ and they are exchanged by the action of $\Z_2^\tau$. Notice that we can identify at which point $e_0$ of the MacLaurin family these two solutions bifurcate since at that point we must have
    $$(\xi_L,\xi_R)=\frac{\Omega}{2}(\mathbf{e_3},-\mathbf{e_3})+\eta_2(e_0)(\mathbf{e_3}, \mathbf{e_3})$$
    and for one of the two components to vanish we must have 
    $$\Omega^2=4 \eta_2^2(e_0)$$
    which is solved for $e_0\simeq 0.8126700$, the right eccentricity value for which the Jacobi and Dedekind's solutions are known to bifurcate (see for instance 33.(b) and chapter 6 in \cite{Cha87}).
\end{itemize}

\paragraph{Acknowledgments:} This work has been partially supported by the Spanish Ministry of Science and Innovation, grant PID2021-125515NB-C21.

\end{document}